\documentstyle[12pt]{article}
\newcommand{\const}{\mathop{\rm const}\limits}

\newcommand{\supp}{\mathop{\rm supp}\limits}

\newcommand{\diam}{\mathop{\rm diam}\limits}

\newcommand{\mes}{\mathop{\rm mes}\limits}

\newcommand{\cov}{\mathop{\rm cov}\limits}

\begin{document}

\begin{center}

{\bf  SIMPLIFICATION OF THE  MAJORIZING MEASURES METHOD,} \par

\vspace{4mm}

{\bf with development.}\par

\vspace{4mm}

 $ {\bf E.Ostrovsky^a, \ \ L.Sirota^b } $ \\

\vspace{4mm}

$ ^a $ Corresponding Author. Department of Mathematics and computer science, Bar-Ilan University, 84105, Ramat Gan, Israel.\\
\end{center}
E - mail: \ galo@list.ru \  eugostrovsky@list.ru\\
\begin{center}
$ ^b $  Department of Mathematics and computer science. Bar-Ilan University,
84105, Ramat Gan, Israel.\\

E - mail: \ sirota3@bezeqint.net\\

\vspace{3mm}
                    {\sc Abstract.}\\

 \end{center}

 \vspace{4mm}

 We update, specify, review and develop in this article the classical majorizing measures method
for investigation of the local structure of random fields,
belonging to X.Fernique and M.Talagrand in order to simplify  and improve the constant values.\par
 Our considerations based on the generalization of the L.Arnold and P.Imkeller generalization of classical
Garsia - Rodemich - Rumsey inequality.
  \vspace{4mm}

\par

\vspace{4mm}

 \par

\vspace{4mm}

\section{Notations. Statement of problem.}

\vspace{3mm}

 Let $ (X,d),  (Y,\rho) $ be separable metric spaces, $ m  $ be arbitrary distribution, i.e.
 Radon probabilistic measure on the  set $ X, $
$ f: X \to Y $  be (measurable) function. Let also $ \Phi(z), \ z \ge 0 $ be continuous Young-Orlicz function,  i.e.
strictly increasing function such that

$$
\Phi(z) = 0 \ \Leftrightarrow z = 0; \ \lim_{z \to \infty} \Phi(z) = \infty.
$$
 We denote as usually

$$
\Phi^{-1}(w) = \sup \{z, z \ge 0, \  \Phi(z) \le w \}, \ w \ge 0
$$
the inverse function to the function $  \Phi; $
$$
B(r,x) = \{ x_1: \ x_1 \in X, \ d(x_1,x) \le r  \}, \ x \in X, \ 0 \le r \le \diam(X)
$$
be the closed ball of radii $ r $ with center at the point $ x. $\par
 Let us introduce the Orlicz space $  L(\Phi) = L(\Phi; m \times m, \ X \otimes X)  $ on the set $ X \otimes X $ equipped with
the Young - Orlicz function $ \Phi. $ \par

 \vspace{4mm}

{\it We assume henceforth that for all the values } $ x_1, x_2 \in X, \ x_1 \ne x_2 $ {\it (the case } $ x_1 = x_2 $  {\it is trivial)}
{\it the value  } $ \rho(f(x_1), f(x_2)) $ {\it belongs to the space } $ L(\Phi). $ \par

\vspace{4mm}
 Note that for the existence of such a function $ \Phi(\cdot) $ is necessary and sufficient only the integrability
of the distance   $ \rho(f(x_1), f(x_2)) $ over the product measure $ m \times m:$

$$
\int_X \int_X \rho(f(x_1), f(x_2)) \ m(dx_1) \  m(dx_2) < \infty,
$$
see \cite{Krasnoselsky1}, chapter 2, section 8. \par

  Under this assumption the distance $ d = d(x_1, x_2) $ may be constructively defined by the formula:

$$
d_{\Phi}(x_1,x_2) := || \rho(f(x_1), f(x_2))||L(\Phi). \eqno(1.1)
$$

  Since the function $ \Phi = \Phi(z) $ is presumed to be continuous and strictly  increasing,
it follows from the relation (1.1) that $ V(d_{\Phi}) \le 1, $ where by definition

$$
V(d):= \int_X \int_X \Phi \left[ \frac{\rho(f(x_1), f(x_2))}{d(x_1,x_2)} \right] \ m(dx_1) \  m(dx_2). \eqno(1.2)
$$

 Let us define also the following important distance function: $ w(x_1, x_2) = $
 $$
  w(x_1, x_2; V) = w(x_1, x_2; V, m ) = w(x_1, x_2; V, m,\Phi) = w(x_1, x_2; V, m,\Phi,d) \stackrel{def}{=}
 $$

$$
 6 \int_0^{d(x_1, x_2)} \left\{ \Phi^{-1} \left[ \frac{4V}{m^2(B(r,x_1))} \right] +
\Phi^{-1}  \left[ \frac{4V}{m^2(B(r,x_2))} \right] \right\} \ dr, \eqno(1.3)
$$
 where $ m(\cdot) $ is probabilistic Borelian measure on the set $ X.$ \par
  The triangle inequality and other properties of the distance function $ w = w(x_1, x_2) $ are proved in
\cite{Kwapien1}.\par

\vspace{3mm}

{\bf Definition 1.1. } (See  \cite{Kwapien1}). The measure $ m  $ is said to be
{\it minorizing measure } relative the distance $ d = d(x_1,x_2), $ if for each values $ x_1, x_2 \in X
\ V(d) < \infty $ and moreover $ \ w(x_1,x_2; V(d)) < \infty. $\par

\vspace{3mm}

We will denote the set of all minorizing measures on the metric set $ (X,d) $  by  $ \cal{M} = \cal{M}(X).$ \par

 Evidently, if the function $ w(x_1, x_2) $ is bounded, then the minorizing measure $  m  $ is majorizing.  Inverse
proposition is not true, see  \cite{Kwapien1}, \cite{Arnold1}. \par

\vspace{3mm}

{\bf Remark 1.1.} If the measure $  m  $ is minorizing, then

$$
w(x_n, x ; V(d)) \to 0 \ \Leftrightarrow d(x_n, x) \to 0, \ n \to \infty.
$$
 Therefore, the continuity of a function relative the distance $  d  $ is equivalent to
the continuity  of this function  relative the distance $  w.  $ \par

\vspace{3mm}

{\bf Remark 1.2.}  If

$$
\sup_{x_1, x_2 \in X} w(x_1, x_2; V(d)) < \infty,
$$
then the measure $ m $ is called {\it majorizing measure.} This classical definition
with theory explanation and applications basically in the investigation of local structure
of random processes and fields  belongs to
X.Fernique \cite{Fernique1},  \cite{Fernique2},  \cite{Fernique3} and M.Talagrand
\cite{Talagrand1}, \cite{Talagrand2}, \cite{Talagrand3}, \cite{Talagrand4}, \cite{Talagrand5}.
 See also \cite{Bednorz1}, \cite{Bednorz2}, \cite{Bednorz3}, \cite{Dudley1}, \cite{Ledoux1},
 \cite{Ostrovsky100}, \cite{Ostrovsky101}, \cite{Ostrovsky102}. \par

\vspace{4mm}

 The following important inequality belongs to  L.Arnold and P.Imkeller \cite{Arnold1}, \cite{Imkeller1};
see also \cite{Kassman1}, \cite{Barlow1}. \par

\vspace{4mm}

{\bf  Theorem of  L.Arnold and P.Imkeller.} {\it Let the measure $  m  $ be minorizing. Then there exists a
modification of the function $ f $ on the set of zero measure, which we denote also by $ f, $ for which}

$$
\rho(f(x_1), f(x_2)) \le  w(x_1, x_2; V, m,\Phi,d). \eqno(1.4)
$$
{\it  As a consequence: this function $  f  $ is $  d - $ continuous and moreover $ w - $ Lipshitz continuous
with unit constant. }\par

\vspace{4mm}

 The inequality (1.4) of L.Arnold and P.Imkeller is significant generalization of celebrated
Garsia - Rodemich - Rumsey inequality, see \cite{Garsia1}, with at the same applications as mentioned before
\cite{Hu1}, \cite{Ostrovsky100}, \cite{Ostrovsky101}, \cite{Ostrovsky102}, \cite{Ral'chenko1}. \par

\vspace{4mm}

{\bf   Our purpose in this report is application of the  L.Arnold and P.Imkeller inequality
to the investigation of continuity of random fields; limit theorems, in particular, Central Limit Theorem,
for random processes; exponential estimates for distribution of maximum of random fields etc. } \par

\vspace{3mm}

 Obtained results improve and generalize recent ones in \cite{Bednorz2}, \cite{Dudley1}, \cite{Fernique1}, \cite{Hu1},
\cite{Ostrovsky1}, \cite{Ostrovsky103},\cite{Ral'chenko1},  \cite{Talagrand2}, \cite{Talagrand4} etc. \par

\vspace{3mm}

{\bf Remark 1.3.} The inequality of  L.Arnold and P.Imkeller (1.4) is closely related with the theory of
fractional order Sobolev's - rearrangement invariant spaces, see  \cite{Barlow1}, \cite{Garsia1}, \cite{Hu1}, \cite{Kassman1},
\cite{Nezzaa1}, \cite{Ostrovsky101}, \cite{Ral'chenko1}, \cite{Runst1}.\par

\vspace{3mm}

{\bf Remark 1.4.} In the previous articles \cite{Kwapien1}, \cite{Bednorz4} was imposed on the function $ \Phi(\cdot) $
the following $ \Delta^2 $ condition:

$$
\Phi(x) \Phi(y) \le \Phi(K(x+y)), \ \exists K = \const \in (1,\infty), \ x,y \ge 0
$$
or equally

$$
\sup_{x,y > 0} \left[ \frac{\Phi^{-1}(xy)}{\Phi^{-1}(x) + \Phi^{-1}(y)}\right] < \infty. \eqno(1.5)
$$
{\it  We do not suppose this condition. For instance, we can consider the function of a view $ \Phi(z) = |z|^p,  $
which does not satisfy (1.5)}. \par

\vspace{3mm}

{\bf Remark 1.5.}  In the works of  M.Ledoux and  M.Talagrand \cite{Ledoux1}, \cite{Talagrand1} -  \cite{Talagrand3}
was investigated at first the case when the function $ d = d(x_1,x_2)  $ is {\it ultrametric }, i.e. satisfies the condition

$$
d(x_1,x_3) \le \max(d(x_1,x_2), d(x_2,x_3) ).
$$
 Note that in the classical monograph of N.Bourbaki \cite{Bourbaki1}, chapter 3, section 5  these function are called pseudometric. \par
We do not use this approach.\par

\vspace{3mm}
\section{ Main result: Continuity of random fields. }

\vspace{3mm}

\begin{center}
{\bf General Orlicz approach.}\\
\end{center}

\vspace{4mm}

  Let $ \xi = \xi(x), \ x \in X $ be separable continuous {\it in probability }
  random field (r.f),  not necessary to be Gaussian.
   The correspondent  set of elementary events,
   probability and expectation we will denote  by $ \omega, \ {\bf P}, \ {\bf E,}  $  and  the probabilistic
 Lebesgue-Riesz $ L_p $   norm of a random variable (r.v) $  \eta $  we will denote as follows:

 $$
 |\eta|_p \stackrel{def}{=} \left[ {\bf E} |\eta|^p \right]^{1/p}.
 $$

   We find in this section some sufficient condition for continuity of $ \xi(x) $  and estimates for its  modulus of continuity
 $ \omega(f,\delta): $

 $$
 \omega(f, \delta) =  \omega(f, \delta,d) := \sup_{x_1,x_2 \in X, d(x_1,x_2) \le \delta}  w(f(x_1), f(x_2)). \eqno(2.0)
 $$

    Recall that the first publication about fractional  Sobolev's  inequalities
\cite{Garsia1}  was devoted in particular to the such a problem; see also articles \cite{Hu1}, \cite{Ostrovsky101}, \cite{Ral'chenko1}. \par
\vspace{3mm}

  Let $ \Phi = \Phi(u) $ be again the Young-Orlicz function. We will denote the Orlicz norm by means of the function $  \Phi $
 of a r.v. $ \kappa $ defined on our probabilistic space $ (\Omega, {\bf P}) $ as $ |||\kappa|||L(\Omega,\Phi) $   or for simplicity
 $ |||\kappa|||\Phi. $  \par
  We introduce  the so-called {\it natural} distance $ \rho_{\Phi}(x_1,x_2) $ as follows:

  $$
  d := d_{\Phi} = d_{\Phi}(x_1,x_2):= |||\rho(\xi(x_1),\xi(x_2))|||L(\Omega,\Phi), \ x_1,x_2 \in X.  \eqno(2.1)
  $$

\vspace{4mm}

{\bf Theorem 2.1.} {\it  Let  $ m(\cdot) $ be the probabilistic minorizing measure on the set $  X  $
relative the distance $  d_{\Phi}(\cdot, \cdot). $ There exists a non-negative random variable $  Z = Z(d_{\Phi},m) $
with unit expectation: $ {\bf E} Z = 1 $  for which }

$$
\rho(\xi(x_1), \xi(x_2)) \le  w(x_1, x_2; Z(d_{\Phi},m)). \eqno(2.2)
$$
  {\it As  a consequence:  the r.f.  $  \xi = \xi(x)  $ is $ d \ - $ continuous with probability one. } \par
\vspace{4mm}

{\bf Proof.}  We pick

$$
Z = \int_X \int_X \Phi \left( \frac{\rho(\xi(x_1), \xi(x_2))}{d_{\Phi}(x_1,x_2)}  \right) \ m(dx_1) \ m(dx_2).
$$
We have by means of theorem Fatou - Tonelli

$$
 {\bf E} Z  = {\bf E} \int_X \int_X  \Phi \left( \frac{\rho(\xi(x_1),\xi(x_2))}{d_{\Phi}(x_1,x_2)} \right) \ m(dx_1) \ m(dx_2) =
$$

$$
 \int_X \int_X {\bf E} \Phi \left( \frac{\rho(\xi(x_1),\xi(x_2))}{d_{\Phi}(x_1,x_2)} \right) \ m(dx_1) \ m(dx_2) = 1, \eqno(2.3)
$$
since $ \int_X \int_X  m(dx_1) \  m(dx_2) = 1.  $\par
 It remains to apply the L.Arnold and P.Imkeller inequality. \par

\vspace{4mm}

{\bf Corollary 2.1.} {\it Under at the same conditions as in theorem 2.1 }

\vspace{4mm}

$$
\rho(\xi(x_1), \xi(x_2)) \le  \inf_{m \in \cal{M} } w(x_1, x_2; Z(d_{\Phi},m)). \eqno(2.4)
$$
{\it where the non-negative r.v. $ Z(d_{\Phi}),m) $ has unit expectation.} \par

\vspace{4mm}

\begin{center}
{\bf Examples.}\par
\end{center}

\vspace{3mm}

{\bf Example 2.1. Lebesgue-Riesz spaces approach.} \par

\vspace{3mm}

 Suppose  the measure $ m $ and distance $  d  $ are such that

$$
|\rho(\xi(x_1), \xi(x_2))|_p \le d(x_1,x_2), \ p = \const \ge 1, \eqno(2.5)
$$

$$
 m^2(B(r,x)) \ge r^{\theta}/C(\theta), \ r \in [0,1], \ \theta = \const > 0, \ C(\theta) \in (0,\infty). \eqno(2.6)
$$

 Let also $ p = \const  > \theta. $\par

 \vspace{3mm}

 {\bf Proposition 2.1.}
{\it We get using the inference of theorem 2.1 that for the r.f. $ \xi = \xi(x) $ the following inequality holds:}
$ m \in  \cal{M} $  and

$$
\rho(\xi(x_1), \xi(x_2) ) \le 12 \ Z^{1/p} \ 4^{1/p} \ C^{1/p}(\theta) \  \frac{d^{1-\theta/p}(x_1, x_2)}{1-\theta/p}, \eqno(2.7)
$$
where the r.v. $  Z  $ has unit expectation:  $ {\bf E} Z = 1. $\par

\vspace{3mm}

{\bf Example 2.2. Grand Lebesgue spaces approach.} \par

\vspace{3mm}

 We recall first of all briefly  the definition ans some simple properties of the so-called Grand Lebesgue spaces;   more detail
investigation of these spaces see in \cite{Fiorenza3}, \cite{Iwaniec2}, \cite{Kozachenko1}, \cite{Liflyand1}, \cite{Ostrovsky1},
\cite{Ostrovsky2}; see also reference therein.\par

  Recently  appear the so-called Grand Lebesgue Spaces $ GLS = G(\psi) =G\psi =
 G(\psi; A,B), \ A,B = \const, A \ge 1, A < B \le \infty, $ spaces consisting
 on all  the random variables (measurable functions) $ f: \Omega \to R $ with finite norms

$$
   ||f||G(\psi) \stackrel{def}{=} \sup_{p \in (A,B)} \left[ |f|_p /\psi(p) \right].
$$

  Here $ \psi(\cdot) $ is some continuous positive on the {\it open} interval
$ (A,B) $ function such that

$$
     \inf_{p \in (A,B)} \psi(p) > 0, \ \psi(p) = \infty, \ p \notin (A,B).
$$
 We will denote
$$
 \supp (\psi) \stackrel{def}{=} (A,B) = \{p: \psi(p) < \infty, \}
$$

The set of all $ \psi $  functions with support $ \supp (\psi)= (A,B) $ will be
denoted by $ \Psi(A,B). $ \par
  This spaces are rearrangement invariant, see \cite{Bennet1}, and
    are used, for example, in the theory of probability  \cite{Kozachenko1},
  \cite{Ostrovsky1}, \cite{Ostrovsky2}; theory of Partial Differential Equations \cite{Fiorenza3},
  \cite{Iwaniec2};  functional analysis \cite{Fiorenza3}, \cite{Iwaniec2},  \cite{Liflyand1},
  \cite{Ostrovsky2}; theory of Fourier series, theory of martingales, mathematical statistics,
  theory of approximation etc.\par

  Notice that in the case when $ \psi(\cdot) \in \Psi(A,\infty)  $ and a function
 $ p \to p \cdot \log \psi(p) $ is convex,  then the space
$ G\psi $ coincides with some {\it exponential} Orlicz space. \par
 Conversely, if $ B < \infty, $ then the space $ G\psi(A,B) $ does  not coincides with
 the classical rearrangement invariant spaces: Orlicz, Lorentz, Marcinkiewicz  etc.\par

  The fundamental function of these spaces $ \phi(G(\psi), \delta) = ||I_A ||G(\psi), \mes(A) = \delta, \ \delta > 0, $
where  $ I_A  $ denotes as ordinary the indicator function of the measurable set $ A, $ by the formulae

$$
\phi(G(\psi), \delta) = \sup_{ p \in \supp (\psi)} \left[ \frac{\delta^{1/p}}{\psi(p)} \right].  \eqno(2.8)
$$
The fundamental function of arbitrary rearrangement invariant spaces plays very important role in functional analysis,
theory of Fourier series and transform \cite{Bennet1} as well as in our further narration. \par

 Many examples of fundamental functions for some $ G\psi $ spaces are calculated in  \cite{Ostrovsky1}, \cite{Ostrovsky2}.\par

\vspace{3mm}

{\bf Remark 2.1} If we introduce the {\it discontinuous} function

$$
\psi_{(r)}(p) = 1, \ p = r; \psi_{(r)}(p) = \infty, \ p \ne r, \ p,r \in (A,B)
$$
and define formally  $ C/\infty = 0, \ C = \const \in R^1, $ then  the norm
in the space $ G(\psi_r) $ coincides with the $ L_r $ norm:

$$
||f||G(\psi_{(r)}) = |f|_r.
$$
Thus, the Grand Lebesgue Spaces are direct generalization of the
classical exponential Orlicz's spaces and Lebesgue spaces $ L_r. $ \par

\vspace{3mm}

{\bf Remark 2.2}  The function $ \psi(\cdot) $ may be generated as follows. Let $ \xi = \xi(x)$
be some measurable function: $ \xi: X \to R $ such that $ \exists  (A,B):
1 \le A < B \le \infty, \ \forall p \in (A,B) \ |\xi|_p < \infty. $ Then we can
choose

$$
\psi(p) = \psi_{\xi}(p) = |\xi|_p.
$$

 Analogously let $ \xi(t,\cdot) = \xi(t,x), t \in T, \ T $ is arbitrary set,
be some {\it family } $ F = \{ \xi(t, \cdot) \} $ of the measurable functions:
$ \forall t \in T  \ \xi(t,\cdot): X \to R $ such that
$$
 \exists  (A,B): 1 \le A < B \le \infty, \ \sup_{t \in T} \
|\xi(t, \cdot)|_p < \infty.
$$
 Then we can choose

$$
\psi(p) = \psi_{F}(p) = \sup_{t \in T}|\xi(t,\cdot)|_p.
$$
The function $ \psi_F(p) $ may be called as a {\it natural function} for the family $ F. $
This method was used in the probability theory, more exactly, in
the theory of random fields, see \cite{Kozachenko1},\cite{Ostrovsky1}, chapters 3,4. \par

 For instance, the function $ \Phi(\cdot) $ in (2.1) may be introduced by a natural way based on the
family

$$
F_{d,X} = \{ \rho(\xi(x_1),\xi(x_2)) \}, \ x_1,x_2 \in X.
$$

\vspace{3mm}

{\bf Remark 2.3} Note that the so-called {\it exponential} Orlicz spaces are particular cases of
Grand Lebesgue spaces  \cite{Kozachenko1}, \cite{Ostrovsky1}, p. 34-37.  In detail, let the $ N- $
Young-Orlicz function has a view

$$
 N(u) = e^{\mu(u)},
$$
where the function $ u \to \mu(u) $ is convex even twice differentiable function such that

$$
\lim_{u \to \infty} \mu'(u) = \infty.
$$
 Introduce a new function
$$
\psi_{\{N\}}(x) = \exp \left\{ \frac{\left[\log N(e^x) \right]^*}{x}   \right\},
$$
where $  g^*(\cdot) $ denotes the Young-Fenchel transform of the function $  g:  $

$$
g^*(x) = \sup_y (xy - g(y)).
$$
 Conversely,  the $  N  - $ function may be calculated up to equivalence
  through corresponding function $ \psi(\cdot) $  as follows:

 $$
 N(u) = e^{\tilde{\psi}^*(\log |u|) }, \ |u| > 3; \ N(u) = C u^2, |u| \le 3; \  \tilde{\psi}(p) = p \log \psi(p).
 $$
 The Orlicz's space $ L(N) $ over our probabilistic space is equivalent up to sublinear norms equality with
Grand Lebesgue space $ G\psi_{\{N\}}. $ \par

\vspace{3mm}

{\bf Remark 2.4.} The theory of probabilistic {\it exponential} Grand Lebesgue spaces
or equally exponential Orlicz spaces gives a
  very convenient apparatus for investigation of
the r.v. with exponential decreasing tails of distributions. Namely, the non-zero  r.v. $ \eta $ belongs to the
Orlicz space $ L(N), $  where $ N = N(u) $ is function described in equality (1.8), if and only if

$$
{\bf P} (\max(\eta, -\eta) > z) \le \exp(-\mu(C z)), \ z > 1,  \ C = C(N(\cdot), ||\eta||L(N)) \in (0,\infty).
$$
(Orlicz's version). \par
 Analogously may be written a Grand Lebesgue version of this inequality.
 In detail,  if $ 0 < ||\eta||G\psi< \infty,  $ then

 $$
{\bf P} (\max(\eta, -\eta) > z) \le  2\exp \left(- \tilde{\psi}(\log [z /||\eta||G\psi] ) \right), z \ge ||\eta||G\psi.
 $$
  Conversely, if

 $$
{\bf P} (\max(\eta, -\eta) > z) \le  2\exp \left(- \tilde{\psi}(\log [z /K] ) \right), z \ge K,
$$
then $ ||\eta||G\psi \le C(\psi)  \cdot K, \ C(\psi) \in (0,\infty). $  \par

\vspace{3mm}

 A very important subclass of the $ G\psi $ spaces form the so-called $ B(\phi) $ spaces. \par

Let $ \phi = \phi(\lambda), \lambda \in (-\lambda_0, \lambda_0), \ \lambda_0 = \const \in (0, \infty] $ be some even strong
convex which takes positive values for positive arguments twice continuous differentiable function, such that
$$
 \phi(0) = 0, \ \phi^{//}(0) \in(0,\infty), \ \lim_{\lambda \to \lambda_0} \phi(\lambda)/\lambda = \infty. \eqno(2.9)
$$
 We denote the set of all these function as $ \Phi; \ \Phi =
\{ \phi(\cdot) \}. $ \par
 We say that the {\it centered} random variable (r.v) $ \xi = \xi(\omega) $
belongs to the space $ B(\phi), $ if there exists some non-negative constant
$ \tau \ge 0 $ such that

$$
\forall \lambda \in (-\lambda_0, \lambda_0) \ \Rightarrow
{\bf E} \exp(\lambda \xi) \le \exp[ \phi(\lambda \ \tau) ]. \eqno(2.10).
$$
 The minimal value $ \tau $ satisfying (2.10) is called a $ B(\phi) \ $ norm
of the variable $ \xi, $ write
 $$
 ||\xi||B(\phi) = \inf \{ \tau, \ \tau > 0: \ \forall \lambda \ \Rightarrow
 {\bf E}\exp(\lambda \xi) \le \exp(\phi(\lambda \ \tau)) \}.
 $$
 This spaces are very convenient for the investigation of the r.v. having a
exponential decreasing tail of distribution, for instance, for investigation of the limit theorem,
the exponential bounds of distribution for sums of random variables,
non-asymptotical properties, problem of continuous of random fields,
study of Central Limit Theorem in the Banach space etc.\par

  The space $ B(\phi) $ with respect to the norm $ || \cdot ||B(\phi) $ and
ordinary operations is a Banach space which is isomorphic to the subspace
consisted on all the centered variables of Orlicz's space $ (\Omega,F,{\bf P}), N(\cdot) $ with $ N \ - $ function

$$
N(u) = \exp(\phi^*(u)) - 1, \ \phi^*(u) = \sup_{\lambda} (\lambda u -\phi(\lambda)).
$$
 The transform $ \phi \to \phi^* $ is called Young-Fenchel transform. The proof of considered assertion used the
properties of saddle-point method and theorem of Fenchel-Moraux:
$$
\phi^{**} = \phi.
$$

 The next facts about the $ B(\phi) $ spaces are proved in \cite{Kozachenko1}, \cite{Ostrovsky1}, p. 19-40:

$$
{\bf 1.} \ \xi \in B(\phi) \Leftrightarrow {\bf E } \xi = 0, \ {\bf and} \ \exists C = \const > 0,
$$

$$
U(\xi,x) \le \exp(-\phi^*(Cx)), x \ge 0,
$$
where $ U(\xi,x)$ denotes in this article the {\it tail} of
distribution of the r.v. $ \xi: $

$$
U(\xi,x) = \max \left( {\bf P}(\xi > x), \ {\bf P}(\xi < - x) \right),
\ x \ge 0,
$$
and this estimation is in general case asymptotically exact. \par
 Here and further $ C, C_j, C(i) $ will denote the non-essentially positive
finite "constructive" constants.\par

 The function $ \phi(\cdot) $ may be "constructively" introduced by the formula
$$
\phi(\lambda) = \phi_0(\lambda) \stackrel{def}{=} \log \sup_{t \in T}
 {\bf E} \exp(\lambda \xi(t)), \eqno(2.11)
$$
 if obviously the family of the centered r.v. $ \{ \xi(t), \ t \in T \} $ satisfies the {\it uniform } Kramer's condition:
$$
\exists \mu \in (0, \infty), \ \sup_{t \in T} U(\xi(t), \ x) \le \exp(-\mu \ x),
\ x \ge 0.
$$
 In this case, i.e. in the case the choice the function $ \phi(\cdot) $ by the
formula (2.11), we will call the function $ \phi(\lambda) = \phi_0(\lambda) $
a {\it natural } function. \par
 {\bf 2.} We define $ \psi(p) = \psi_{\phi}(p) := p/\phi^{-1}(p), \ p \ge 2. $
   It is proved that the spaces $ B(\phi) $ and $ G(\psi) $ coincides:$ B(\phi) =
G(\psi) $ (set equality) and both
the norm $ ||\cdot||B(\phi) $ and $ ||\cdot|| $ are equivalent: $ \exists C_1 =
C_1(\phi), C_2 = C_2(\phi) = \const \in (0,\infty), \ \forall \xi \in B(\phi) $

$$
||\xi||G(\psi) \le C_1 \ ||\xi||B(\phi) \le C_2 \ ||\xi||G(\psi).
$$

   The Gaussian (more precisely, subgaussian) case considered in \cite{Garsia1}, \cite{Hu1}, \cite{Ral'chenko1}
 may be obtained by choosing $ \Phi(z) = \Phi_2(z) := \exp(z^2/2) - 1 $ or equally $ \psi(p) = \psi_2(p) = \sqrt{p}. $
 It may be considered easily more general example when
 $ \Phi(z) =  \Phi_Q(z) := \exp(|z|^Q/Q) - 1, \ Q = \const > 0; \ \Leftrightarrow \psi(p) = \psi_Q(p) := p^{1/Q}, \ p \ge 1. $  \par

 In the last case the following implication holds:

 $$
 \eta \in L(\Phi_Q), \ Q > 1 \ \Leftrightarrow  U(\eta,x) \le \exp \left( - C(\Phi, \eta) \ x^{Q'}  \right),
 $$
where as usually $ Q' = Q/(Q-1). $\par

\vspace{4mm}

 Assume that  the number $  \theta, $ measure $ m, $  distance $  d,  $ and the function $ \psi = \psi(p) $  are such that
$  \theta > 0, \ C(\theta) = \const \in (0,\infty); $

$$
(A,B):= \supp \psi(\cdot), \ \tilde{A}:= \max(A,\theta), \ B > \tilde{A}; \eqno(2.12)
$$

$$
||\rho(\xi(x_1), \xi(x_2))||G\psi \le d(x_1,x_2); \eqno(2.13)
$$

$$
 m^2(B(r,x)) \ge r^{\theta}/C(\theta), \ r \in [0,1], \  \ C(\theta) \in (0,\infty). \eqno(2.14)
$$

 Define also a new   function:

 $$
\psi_{\theta}(p) \stackrel{def}{=} (1-\theta/p) \ \psi(p), \ p \in (\tilde{A}, B). \eqno(2.15)
 $$

 \vspace{3mm}

 {\bf Proposition 2.2.} {\it  Under formulated above conditions (2.12), (2.13), (2.14) we have:}
  $ m \in  \cal{M} $  {\it and}

$$
||\rho(\xi(x_1), \xi(x_2) )||G\psi \le 12 \ d(x_1,x_2) \   \phi \left(G\psi_{\theta},  4 \ C(\theta) \ d^{-\theta}(x_1,x_2) \right). \eqno(2.16)
$$
{\bf Proof.}  The condition (2.13) implies that

$$
|\rho(\xi(x_1), \xi(x_2))|_p \le \psi(p) \cdot d(x_1,x_2), \ p \in (\tilde{A}, B).
$$
 We derive using proposition 2.1

 $$
\frac{ |\rho(\xi(x_1), \xi(x_2) )|_p}{12 \ \psi(p) \ d(x_1,x_2) } \le \frac{(4  C(\theta))^{1/p} \ d^{-\theta/p}(x_1,x_2)}{(1-\theta/p) \psi(p)} =
 \frac{(4  C(\theta))^{1/p} \ d^{-\theta/p}(x_1,x_2)}{ \psi_{\theta}(p)}. \eqno(2.17)
$$
 The assertion (2.16) of proposition 2.2 follows immediately on the basis of the definition of fundamental function for Grand
Lebesgue spaces  from (2.17) after taking supremum over $ p. $\par

\vspace{3mm}

{\bf Remark 2.5.}  The case when $ \psi(p) = \sqrt{p}  $ appropriates to the Gaussian (more generally, subgaussian) random field $ \xi(x). $
The case $ \psi(p) = \exp(C p) $ appears in the articles \cite{Arnold1} and  \cite{Imkeller1}. However, in both these cases the
condition (1.5) is satisfied.\par
 In the case $ \psi(p) = \psi_{(r)}(p) $ we obtain the proposition 2.1. as a  particular case. \par

\vspace{4mm}

  Obtained in this section results specify and generalize  ones  in the articles \cite{Garsia1}, \cite{Hu1}, \cite{Ral'chenko1}. \par
 Another approach to the problem of (ordinary) continuity  of random fields based on the so-called generic chaining method
 and entropy  technique with described applications see in  \cite{Bednorz1}, \cite{Fernique1}, \cite{Kozachenko1}, \cite{Ledoux1},
 \cite{Ostrovsky1}, \cite{Talagrand1}, \cite{Talagrand2} etc.\par

\vspace{3mm}

\section{Weak compactness of random fields.}

\vspace{3mm}

\begin{center}

{\bf A. General result. }\\

\end{center}

\vspace{3mm}

 Let $  \xi_n = \xi_n(x), \ x \in X, \ n = 1,2,\ldots $ be a {\it  sequence } of separable stochastic continuous random fields.  \par
  {\it  We suppose in this section that the metric space $ (X,d) $ is compact and  that there exists
   a non-random point $ x_0 \in X $  for which the one-dimensional r.v.  $ \xi_n(x_0) $ are tight.} \par
 This condition is satisfied  if for instance $  Y = R^1, \ {\bf E} \xi_n(x_0) =  0  $ and $ \sup_n {\bf Var} \ \xi_n(x_0) < \infty. $\par

\vspace{3mm}

 Let $ \overline{\Phi} = \overline{\Phi}(u) $ be again the Young-Orlicz function, a single function for all the r.f. $ \xi_n(x). $
Define the common distance

$$
\overline{d}(x_1,x_2) =  \sup_n d_n(x_1,x_2),   \ x_1,x_2 \in X,
$$
where
$$
   d_n(x_1,x_2)= ||\rho( \xi_n(x_1),\xi_n(x_2)) || L(\overline{\Phi}).\eqno(3.1)
$$

  Suppose that the set $  X  $ is compact set relative the distance $ \overline{d} = \overline{d}(x_1,x_2). $ \par

  \vspace{3mm}

 We intend in this section to obtain the sufficient condition for weak compactness of the distributions of the
family  $  \xi_n = \xi_n(\cdot) $ in the space of continuous functions $  C(X) = C(X; \overline{d}). $\par

\vspace{4mm}

{\bf Theorem 3.1.} {\it  Let in addition $ m(\cdot) $ be probabilistic minorizing measure on the set $  X  $
relative the distance } $  \overline{d}(\cdot, \cdot) $  {\it and  the Young - Orlicz function $ \overline{\Phi}. $
  Then the family of distributions generated by
the continuous versions of r.f.  } $  \{  \xi_n(\cdot) \} $ {\it in the space of $ \overline{d} $  continuous functions }
$ C(X) = C(X, \overline{d}): $

$$
\nu_n(A) = {\bf P} ( \xi_n(\cdot) \in A),
$$
{\it where $ A $ is arbitrary Borelian set in $ C(X), $ is weakly compact, i.e. tight.  }\par

\vspace{3mm}
{\bf Proof.} It follows from theorem 2.1 that
there exists a sequence of non-negative random variable $  Y_n = Y_n(d_{\overline{\Phi}},m) $
with unit expectation: $ {\bf E} Y_n = 1 $  for which

$$
\overline{\rho}(\xi_n(x_1), \xi_n(x_2)) \le  w(x_1, x_2; Y_n,\overline{d},m). \eqno(3.2)
$$
 Therefore,

$$
 \forall \epsilon > 0 \ \Rightarrow  \lim_{\delta \to 0+} {\bf P}(\omega(\xi_n,\delta) > \epsilon) = 0. \eqno(3.3)
$$

 This completes the proof of theorem 3.1. \par
\vspace{3mm}

{\bf Consequence 3.1.}  Suppose in addition to the conditions of theorem 3.1 that as $ n \to \infty $ the
finite-dimensional distributions of the r.f. $ \xi_n(x) $ converge to the
finite-dimensional distributions of some  r.f. $ \xi_{\infty}(x).  $ Then the sequence of distributions
$ \nu_n(\cdot) $ weakly converges to the $ \nu_{\infty}(\cdot). $ Namely, for every continuous bounded functional
$ F: C(X,\overline{d}) \to R $

$$
\lim_{n \to \infty}  {\bf E} F(\xi_n(\cdot)) = {\bf E} F(\xi_{\infty}(\cdot)). \eqno(3.4)
$$
\vspace{3mm}

{\bf Remark 3.1.} Let $ \Phi_n(u)  $ be "individual" Young - Orlicz function for each field $  \xi_n(\cdot) $
described below. The "common" Young - Orlicz function  $ \overline{\Phi}(u)  $ may be constructed  evidently as follows:

$$
\overline{\Phi}(u) = \sup_n \Phi_n(u),
$$
if it is finite for all the values $ u, \ u \in R. $ \par
 For instance, let $ \psi_n = \psi_n(p) $ be natural function for the r.f. $ \xi_n(x); \ 1 \le p \le b_n. $
Assume $ b:= \inf_n b_n  > 1 $ and suppose

$$
\psi_{\infty}(p) := \sup_n \psi_n(p) < \infty, \ 1 \le p < b. \eqno(3.5)
$$
 Then the $ G\psi_{\infty} $ space is suitable for us; the correspondent  Young - Orlicz function $ \Phi_{\infty} = \overline{\Phi} $ is described
in the second section. \par
 Analogously  may be considered the case of $ B(\phi)  $ spaces.
 The common function $  \overline{\phi}(\cdot) $ may be "constructively" introduced by the formula
$$
\overline{\phi}(\lambda) = \phi_{\infty}(\lambda) \stackrel{def}{=} \log \sup_n \sup_{x \in X}
 {\bf E} \exp(\lambda \xi_n(x)),  \eqno(3.6)
$$
 if obviously the family of the centered r.f. $ \{ \xi_n(x), \ n= 1,2, \ldots; \ x \in X \} $ satisfies the {\it uniform } Kramer's condition:
$$
\exists \mu \in (0, \infty), \ \sup_n \sup_{x \in X} U(\xi_n(x), \ u) \le \exp(-\mu \ u), \ u \ge 0. \eqno(3.7)
$$

\vspace{3mm}
\begin{center}

{\bf B. Central Limit Theorem in the space of continuous functions. }\\

\end{center}
\vspace{3mm}

  In this subsection $ Y = R^1 $ and $ \rho $  is any continuous distance in $  Y. $ \par

 Let $ \eta_i = \eta_i(x), \ x \in X $ be independent centered: $ {\bf E} \eta_i(x) = 0 $ identical  distributed
random fields with finite covariation function $ R(x_1,x_2) = \cov(\eta_i(x_1), \eta_i(x_2)) = {\bf E}\eta_i(x_1) \cdot \eta_i(x_2).$
 Denote

$$
S_n(x) = n^{-1/2} \sum_{i=1}^n \eta_i(x). \eqno(3.8)
$$
 Obviously, the finite-dimensional distributions of r.f. $ S_n(\cdot) $ converge as $ n \to \infty $
to the finite-dimensional distributions of the centered Gaussian r.f. $ S_{\infty}(\cdot) $  with at the same
covariation function $ R(\cdot, \cdot). $\par
\vspace{3mm}
{\bf Definition 3.1.} (See \cite{Dudley1}, \cite{Kozachenko1}). \\

\vspace{2mm}

 We will say that the sequence of r.f. $ \{ \eta_i(\cdot) \}  $ satisfies the Central Limit
Theorem (CLT) in the space $ C(X,d), $  if  $ {\bf P}\left(\eta_i(\cdot) \in C(X,d) \right)  = 1 $ and
the distributions of the r.f. $ S_n(\cdot)  $ in the set $ C(X,d) $ converge weakly as $ n \to \infty $ to the
distribution of  the Gaussian r.f. $ S_{\infty}(\cdot). $\par
\vspace{3mm}

 We formulate here  some  sufficient conditions  for the CLT in the space of continuous functions in the terms
of minorizing  measures. \par
 Note that in the terms of majorizing measures these conditions are obtained, e.g., in
\cite{Heinkel1}, \cite{Dudley1}, \cite{Ledoux1}, \cite{Talagrand2}; in the  entropy terms - in
\cite{Kozachenko1}, \cite{Ostrovsky1}, chapter 4, section 4 etc. \par

 We can use the the result of the last subsection.  Namely, let $ \zeta = \zeta(\lambda) $ be
natural $ \phi \ - $  function for the r.f. $ \eta_1(x): $

 $$
 \zeta(\lambda) := \sup_{x \in X} \log {\bf E} \exp( \lambda \eta_1(x) ), \eqno(3.9)
 $$
 if as before the family of the centered r.v. $ \{ \eta_1(x), \ x \in X \} $ satisfies the  uniform  Kramer's condition.

 We have:

 $$
\log {\bf E} \exp(\lambda S_n(x))  \le n^{-1/2} \zeta(\lambda/n) \le \overline{\zeta}(\lambda), \eqno(3.10)
 $$
 where by definition

 $$
  \overline{\zeta}(\lambda) \stackrel{def}{=} \sup_n [n^{-1/2} \zeta(\lambda/n) ] < \infty, \ |\lambda| < \lambda_0 = \const > 0. \eqno(3.11)
 $$

 For instance, let

 $$
 \zeta(\lambda) = \lambda^2 \ I(|\lambda| \le 1) +   |\lambda|^Q \ I(|\lambda| > 1), \ Q = \const \ge 1,
 $$
then

$$
\overline{\zeta}(\lambda) \asymp \lambda^2 \ I(|\lambda| \le 1) +   |\lambda|^{\max(Q,2)} \ I(|\lambda| > 1).
$$

 The equivalent conclusion in the terms of $ \Phi - $ functions:  if here

 $$
 \Phi(u) = \exp \left( |u|^{\beta}   \right) - 1, \ |u| \ge 1, \ \beta = \const > 1,
 $$
then

 $$
 \overline{\Phi}(u) \asymp \exp \left( C \ |u|^{\max(\beta',2)'} \right) - 1, \ |u| \ge 1.
 $$

\vspace{4mm}
 Introduce also the Young - Orlicz function

$$
\Theta(u) = e^{[\overline{\zeta}]^*(u)} - 1
$$
and the correspondent distance

$$
\theta(x_1,x_2) = ||\eta_1(x_1) - \eta_1(x_2)|| L(\Theta).
$$

\vspace{4mm}

{\bf Theorem 3.2.} {\it  Let  $ m(\cdot) $ be any probabilistic minorizing measure on the set $  X  $
relative the distance } $  \theta(\cdot, \cdot) $  {\it and  the Young - Orlicz function} $ \Theta. $
{\it Then the sequence of r.f. $ \{ \eta_i(\cdot) \}  $ satisfies the Central Limit
Theorem (CLT) in the space} $ C(X,\theta). $  \par

 \vspace{4mm}
 {\bf Proof} follows immediately from the theorem 3.1., where we set $ \xi_n(x) = S_n(x).  $
 It remains to  ground only the weak compactness of the family r.f. $ S_n(\cdot). $ Note that

 $$
 {\bf E}e^{\lambda S_n(x)} \le e^{\overline{\zeta}(\lambda)},
 $$
or equally

$$
 \sup_n \sup_{x \in X} || S_n(x)|| L(\Theta) < \infty,\eqno(3.12)
$$
see \cite{Ostrovsky1}, chapter 1, section 2.  Analogously,

$$
\sup_n || S_n(x_1) - S_n(x_2)|| L(\Theta) \le C \cdot \theta(x_1, x_2). \eqno(3.13)
$$
 This completes the proof of theorem 3.2; see, e.g. \cite{Prokhorov1}. \par

\vspace{4mm}

 Another approach used the conception of $  G\psi $ spaces.  In detail,  introduce the natural
$  \psi $ function as described below:

$$
\psi(p) = \sup_{x \in X} |\eta_1(x)|_p,
$$
and suppose $ \exists p > 2 $  such that   $ \psi(p) < \infty.  $ \par

 The following distance is finite:

 $$
 d_{\psi}(x_1,x_2) = ||\eta_1(x_1) -\eta_1(x_2)||G\psi.
 $$
 Define a new  $ \psi $ function

$$
\overline{\psi}(p) = \left[\frac{p}{\log(p+1)} \right] \cdot \psi(p).
$$

 It follows from the famous Rosenthal inequality that

 $$
\sup_{x \in X} \sup_n ||S_n||G\overline{\psi}  < \infty
 $$
and

$$
\sup_n ||S_n(x_1) - S_n(x_2) ||G\overline{\psi} \le C \cdot d_{\psi}(x_1,x_2).
$$
 It remains to use the proposition of theorem 3.1. \par

\vspace{4mm}

\section{Non-asymptotical estimates of maximum for random fields.}

\vspace{3mm}

 \begin{center}

 {\bf Grand Lebesgue spaces approach.}

 \end{center}

\vspace{3mm}

 Let  $ \xi = \xi(x), \ x \in X  $ be again separable random field (or process)
with values in the real axis $  R, \ T = \{x \} $ be arbitrary Borelian subset of  $ X.$

 Denote

 $$
 \overline{\xi}_T = \sup_{x \in T} \xi(x), \  \overline{\xi} = \overline{\xi}_X = \sup_{x \in X} \xi(x).
 $$

\vspace{3mm}

{\bf Proposition 4.1.}  Let all the notation of proposition 2.2  be retained and condition   be satisfied.
Denote also

$$
D = \diam(X,d) = \sup_{x_1,x_2 \in X} d(x_1,x_2).
$$
 We assert:

$$
|| \overline{\xi} ||G\psi \le \inf_{x_0 \in X} ||\xi(x_0)||G\psi \ +  \ 12 \ D \cdot \phi \left( G \psi_{\theta}, 4 C(\theta) D^{-\theta}  \right).
\eqno(4.1)
$$
{\bf Proof.} Let $ x_0 $ be arbitrary point in the set $  X; $ we have

$$
\xi(x) = \xi(x_0) + (\xi(x) - \xi(x_0)) \le \xi(x_0) + \rho(\xi(x_0),\xi(x)),
$$

$$
\sup_x \xi(x) \le \xi(x_0) + \sup_x \rho(\xi(x_0),\xi(x)) \le \xi(x_0)  + \sup_{x_1,x_2} \rho(\xi(x_1),\xi(x_2)).
$$
 We conclude using triangle inequality for $ G\psi $ norm  and inequality 2.16:

$$
||\overline{\xi}||G\psi \le ||\xi(x_0)||G\psi + 12 D \cdot \phi \left( G \psi_{\theta}, \ 4 C(\theta) D^{-\theta}  \right).
$$

 Since the value $ x_0 $ is arbitrary, we convince itself that estimate (4.1) is true.\par

{\bf Example 4.1.} Let us consider as a particular case the possibility  $ \psi(p) = \psi_{(r)}(p), $
where $ r = \const > \max(\theta,1). $ We get:

$$
|\overline{\xi}|_r \le \inf_{x_0 \in X} |\xi(x_0)|_r  + 12 \cdot \frac{ (4C(\theta))^{1/r} \cdot D^{1-\theta/r}}{1-\theta/r}.\eqno(4.2)
$$
Here

$$
d(x_1,x_2) = d_r(x_1,x_2)= |\xi(x_1) - \xi(x_2)|_r. \eqno(4.3)
$$

\vspace{3mm}

 \begin{center}

 {\bf Exponential estimations.}

 \end{center}

\vspace{3mm}

  Let $ \Phi = \Phi(u) $ be again as in the beginning of the second section the Young-Orlicz function
generated by the r.f. $  \xi(x), $ so that

$$
\sup_{x \in X} |||\xi(x)|||L(\Phi) < \infty;
$$
in the sequel we will conclude without loss of generality

$$
\sup_{x \in X} ||| \xi(x)|||L(\Phi) = 1. \eqno(4.4)
$$
 Recall that we will denote the Orlicz norm by means of the function $  \Phi $
 of a r.v. $ \kappa $ defined on our probabilistic space  as $ |||\kappa|||L(\Phi).   $ \par
  We introduce  the natural distance $ \rho_{\Phi}(x_1,x_2) $ as follows:

  $$
  d_{\Phi}(x_1,x_2):= ||| \rho(\xi(x_1),\xi(x_2)) |||L(\Phi), \ x_1,x_2 \in X;  \eqno(4.5)
  $$
then in  particular

$$
D = D(\Phi) := \diam(X,d_{\Phi}) \le 2.
 $$

\vspace{4mm}

 For arbitrary Borelian subset $ T \subset X $ we denote
$$
Q(T,u) = {\bf P}( \sup_{t \in T} \xi(t) > u), \ u \ge 2. \eqno(4.6)
$$

$$
Q_+(T,u) = {\bf P}( \sup_{t \in S} |\xi(t)| > u), \ u \ge 2. \eqno(4.7)
$$

  Our purpose in the rest of this section is obtaining an exponentially exact as $ u \to \infty, \ u > u_0 = \const > 0 $
 estimation for the probability $ Q(u) \stackrel{def}{=} Q(X,u), \  Q_+(u) \stackrel{def}{=} Q_+(X,u)  $ in the
 terms of "minorizing measures"  and $ B(\phi) $ spaces.  \par

  In the entropy terms this problem is considered in  \cite{Dudley1}, \cite{Fernique2}, \cite{Fernique3},
 \cite{Ostrovsky1}, chapter 3,  \cite{Ostrovsky103}: in the terms of majorizing measures- in \cite{Ledoux1},
\cite{Ostrovsky102},  \cite{Talagrand1}  etc.

  The estimations of $  Q_+(u) \stackrel{def}{=} Q_+(X,u)  $ are used in the Monte-Carlo method, statistics,
 numerical methods etc., see \cite{Frolov1}, \cite{Grigorjeva1},
 \cite{Ostrovsky1},  \cite{Ostrovsky105}, \cite{Ostrovsky106}. \par

\vspace{3mm}
For instance, we can suppose that the random field $ \xi(x)  $ to be centered and satisfies the uniform Kramer's condition,
so that the natural function

$$
\phi(\lambda) = \log \sup_{x \in X} {\bf E} \exp(\lambda \ \xi(x))
$$
is finite in some non-trivial interval $ \lambda \in (-\lambda_0, \lambda_0), \ \lambda_0 = \const \in (0, \infty].  $\par
 Then we may introduce the following Young-Orlicz function (up to multiplicative positive constant)

 $$
 \Phi_{\phi}(u) = \exp(\phi^*(u))-1,
 $$
so that  $ \sup_{x \in X} ||\xi(x)||B(\phi) = 1  $ and following
$ \sup_x ||\xi(x)||(\Phi) < \infty. $ \par

\vspace{3mm}
 {\it We need also to suppose that the function $ \Phi = \Phi(u) $ satisfies the condition (1.5) with finite
 non-zero constant} $ K. $ \par

\vspace{3mm}

  Let us introduce the following constant (more exactly, functional)

  $$
C_2 = C_2(\Phi)  = \frac{\Phi^{-1}(1)}{54 K^2}, \eqno(4.8)
  $$
and define by $ N(T) = N(T, w,\epsilon), \ \epsilon > 0  $ as  usually for the (pre-compact) metric set $ (T,w), \ T \subset X  $
the minimal number of closed balls with radii $ \epsilon: \  B(x_j, \epsilon) = B(x_j, w,\epsilon) = \{x_1: w(x.x_1) \le \epsilon\}  $
which cover the set $ T: $

$$
T \subset \cup_{j=1}^{N(T)} B(x_j, w,\epsilon); \ N(\epsilon) := N(X,w,\epsilon). \eqno(4.9)
$$

 Recall  that the logarithm  of $ N(X,w,\epsilon) $

 $$
H(X,w,\epsilon)  = \log N(X,w,\epsilon)
 $$
is called "entropy" of the set $  X $ relative the distance $ w(\cdot, \cdot) $ and widely used in the entropy
approach to the investigation of continuity of random processes and fields. \par

\vspace{3mm}

{\bf Proposition 4.2.} {\it  Under formulated above conditions }

$$
Q(u) \le \inf_{\delta \in (0, D)} \frac{N(X, w, \delta)}{\Phi \left( u/(1+\delta/C_2(\Phi))\right)}, \eqno(4.10)
$$

$$
Q_+(u) \le 2 \ \inf_{\delta \in (0, D)} \frac{N(X, w, \delta )}{\Phi \left( u/(1+\delta/C_2(\Phi))\right)}, \eqno(4.11)
$$

\vspace{3mm}

{\bf Proof.} \\
{\bf 1.} S.Kwapien and J.Rosinsky   proved  in \cite{Kwapien1}  the following inequality:

$$
{\bf E}  \Phi \left(2 \ C_2 \sup_{t \ne s} \frac{(\xi(t)-\xi(s))}{w(t,s)} \right) \le 1 + \sup_{t \ne s}
{\bf E} \Phi \left(\frac{(\xi(t) - \xi(s))}{d(t,s)}  \right). \eqno(4.12)
$$

  As long as we choose $ d(t,s) = d_{\Phi}(t,s),  $ we have

$$
{\bf E}  \Phi \left(2 \ C_2 \sup_{t \ne s} \frac{(\xi(t)-\xi(s))}{w(t,s)} \right) \le 2.\eqno(4.13)
$$

 Recall that $ \Phi = \Phi(u) $ is convex function and $ \Phi(0) = 0;  $ following

 $$
 \Phi \left(\frac{u}{2} \right) = \Phi \left(\frac{1}{2} \cdot 0  + \frac{1}{2}\cdot u \right) \le
 \frac{1}{2}  \Phi(0) + \frac{1}{2} \Phi(u) = \frac{1}{2} \Phi(u),
 $$

 We conclude on the basis of inequality (4.13)

$$
{\bf E}  \Phi \left( C_2 \sup_{t \ne s} \frac{(\xi(t)-\xi(s))}{w(t,s)} \right) \le 1,\eqno(4.14)
$$
or equally

$$
||| \ \sup_{w(x_1, x_2) \le \delta} (\xi(x_1) - \xi(x_2)) \ ||| \Phi \le \delta/C_2. \eqno(4.15)
$$

\vspace{3mm}

{\bf 2.} Let $ x_0 $ be fixed (non-random) point in the set $  X. $ Consider the ball
$ B= B(x_0, w, \delta) = \{x_1, x_1 \in X, w(x_0,x_1) \le \delta \}, 0 < \delta \le D. $ We have for the values $ x \in B $
using triangle inequality and (4.15): $ \xi(x)  = \xi(x_0) + [\xi(x) - \xi(x_0)]; $

$$
 \sup_{x \in B}\xi(x) \le \xi(x_0) +
\sup_{w(x_1, x_2) \le \delta} (\xi(x_1) - \xi(x_2)),
$$

$$
||| \ \sup_{x \in B} \xi(x) \ ||| \Phi \le |||\ \xi(x_0) \ |||\Phi + ||| \ \sup_{w(x_1, x_2) \le \delta} (\xi(x_1) - \xi(x_2)) \ ||| \Phi \le
$$

$$
1 + \delta/C_2. \eqno(4.16)
$$
 It follows from Tchebychev's inequality

 $$
 Q(B(x_0, w, \delta) ,u) \le 1/\Phi(u/(1 + \delta/C_2)). \eqno(4.17)
 $$

\vspace{3mm}
{\bf 3.}  The first assertion of proposition 4.2 follows now immediately:

$$
Q(u) = {\bf P} \left[ \cup_{j=1}^{N(T,w,\delta)} \left\{ \sup_{x \in B(x_j,w,\delta)} \xi(x) > u \right\} \right] \le
$$

$$
\sum_{j=1}^{N(T,w,\delta)} {\bf P} \left[ \left\{ \sup_{x \in B(x_j,w,\delta)} \xi(x) > u \right\} \right] \le
N(T,w,\delta) \cdot [ 1/\Phi(u/(1 + \delta/C_2)) ] \eqno(4.18)
$$
after minimization over $ \delta. $\par
 The second assertion of proposition 4.2 follows  from the inequality

 $$
 Q_+(u) \le {\bf P}(\sup_{x \in X}\xi(x) > u) + {\bf P}(\sup_{x \in X}(-\xi(x)) > u).
 $$

\vspace{4mm}

\begin{center}

{\bf Examples.}

\end{center}
\vspace{3mm}
{\bf Example 4.1.} Suppose in addition  to the conditions of proposition 4.2
that the function $ u \to \Phi(u), \ u > 0 $ is logarithmical convex:

$$
(\log \Phi)^{''} (u) > 0.
$$
 Let also $ \gamma = \const \in (0,1).  $ Denote

 $$
 \delta_0 = \delta_0(u; \gamma, \Phi) = \frac{C_2 \gamma}{u \cdot [\log \Phi]'(u)}. \eqno(4.19)
 $$
 We obtain choosing $ \delta = \delta_0 $ substituting into (4.10) that for all sufficiently greatest values
 $ u: \ \delta_0(u; \gamma, \Phi) < D $

$$
Q(u) \le \frac{(1-\gamma)^{-1}}{\Phi(u)} \cdot N \left( \frac{C_2 \gamma}{u \cdot  [\log \Phi]'(u)} \right). \eqno(4.20)
$$

\vspace{4mm}

{\bf Example 4.2.} Let now $ \Phi(u) = \exp(u^2/2) - 1 $ (subgaussian case).  Suppose

$$
N(\epsilon) \le C_3 \epsilon^{-\kappa}, \ \epsilon \in (0,D), \ \kappa = \const > 0.  \eqno(4.21)
$$
 The value $ \kappa $ is said to be {\it majorital} dimension of the set $ X $ relative the distance $ w. $ \par
The optimal value $ \gamma $ in (4.20) if equal to $ \gamma = \gamma_0 := \kappa/(\kappa + 1) $ and we conclude for the
values $ U $ such that

$$
\delta_0 = \frac{C_2 \kappa}{(\kappa + 1) u^2} \le D:
$$

$$
Q(u) \le C_3 \ C_2^{-\kappa} \ \kappa^{-\kappa} \ (\kappa+1)^{ \kappa + 1} \  u^{2 \kappa} \ e^{-u^2/2}. \eqno(4.22)
$$

\vspace{4mm}

{\bf Example 4.3.} Assume that  in the example 4.2 instead the condition (4.21) the following condition holds:

$$
N(\epsilon) \asymp  C_4 e^{ C_5 \ \epsilon^{-\beta}}, \ \epsilon \in (0,D);  \ \beta  = \const > 0.  \eqno(4.23)
$$

 Then

 $$
 Q(u) \le e^{-0.5 u^2 + C_6 u^{2\beta/(\beta+1)} }, \ u \ge C_7. \eqno(4.24)
 $$

 Note that in the  case  $ \beta \ge 2 $ the so-called entropy series

$$
\sum_{n=1}^{\infty} 2^{-n} \ H^{1/2} \left(X,w, 2^{-n} \right)
$$
diverges.\par

\vspace{4mm}

\section{Concluding remarks.}

\vspace{3mm}

{\bf A. Degrees.} \par
 Let $ X = [0,1]^n \ n = 2,3,\ldots. $  In the articles   \cite{Ral'chenko1}, \cite{Hu1}  is obtained a multivariate generalization of
famous Garsia-Rodemich-Rumsey inequality \cite{Garsia1}. Roughly speaking, instead degree "2" in the inequalities
(1.3) and (1.4) stands degree 1  and coefficients dependent on $ d. $ \par
 The ultimate value of this degree in general case of arbitrary metric space $ (X,d) $ is now unknown; see also  \cite{Arnold1},
\cite{Imkeller1}. \par

\vspace{3mm}

{\bf B. Spaces.} \par
 Notice that in all considered cases and under our conditions when
  $ \sup_x ||\xi(x)|| < \infty,  $ then $ ||\sup_x \xi(x)|| < \infty. $ But in the article
\cite{Ostrovsky104} was constructed a "counterexample": there exists a continuous a.e random
process for which

$$
 \sup_x ||\xi(x)|| < \infty, \hspace{5mm}  ||\sup_x \xi(x)|| = \infty.
$$
 This circumstance imply that our conditions are  only sufficient but not necessary.\par

\vspace{3mm}

{\bf C. Rectangle distance.} \par
 In the report \cite{Ostrovsky101} for the multivariate functions  was introduced the rectangle distance.
 For instance, if the r.f. $ \xi = \xi(x,y), \ x,y \in [0,1] $ is bivariate, then

$$
r_{\xi}(x_1,x_2;y_1,y_2) = ||\xi(x_2,y_2) - \xi(x_1,y_2) - \xi(x_1,y_1) + \xi(x_1,x_2)||.
$$
 It is very interest by our opinion to obtain in general case estimations for $ r_{\xi}(x_1,x_2;y_1,y_2)  $
in the terms of minorizing measures alike the one-variate case considered here. \par

\vspace{3mm}

{\bf D.  Lower bounds.} \par
 The lower estimates for probabilities $ Q(u)  $ are obtained e.g. in \cite{Ostrovsky1}, chapter 3, sections 3.5-3.8.
They are  obtained in entropy terms, all the more so in the terms of minorizing measures.  \par
Note that the lower bounds in \cite{Ostrovsky1} may coincide up to multiplicative constants with upper bounds. \par

\vspace{4mm}
{\bf   Acknowledgements.} The authors  would very like to thank  dr. E.Rogover
for Yours remarkable  and  useful remarks and corrections. \par

\end{document}